\documentstyle[12pt]{amsart}

\newcommand{\ie}{{\it i.e.\/}\ }

\numberwithin{equation}{section}
\newtheorem{thm}{Theorem}[section]
\newtheorem{prop}[thm]{Proposition}

\newtheorem{lemma}[thm]{Lemma}

\newcommand{\thmref}[1]{Theorem~\ref{#1}}
\newcommand{\lemref}[1]{Lemma~\ref{#1}}

\begin{document}
\baselineskip=17pt

\title[Parabolic Higgs bundles and Teichm\"uller
spaces]{Parabolic Higgs bundles and Teichm\"uller spaces for
punctured surfaces}
\author[I. Biswas]{Indranil Biswas}
\address{\noindent School of Mathematics, Tata Institute of
Fundamental Research, Bombay, India}
\email{indranil@@math.tifr.res.in}
\author[P. A. Gastesi]{Pablo A. Gastesi}
\address{School of Mathematics, Tata Institute of Fundamental
Research, India}
\curraddr{The Institute of Mathematical Sciences, Madras, India}
\email{pablo@@imsc.ernet.in}
\author[S. Govindarajan]{Suresh Govindarajan}
\address{Theoretical Physics Group, Tata Institute of Fundamental Research,
Bombay, India}
\curraddr{Dept. of Physics, Indian Institute of Technology,
Madras, India}
\email{suresh@@imsc.ernet.in}


\date{\today}

\maketitle
\begin{abstract}
In this paper we study the relation between parabolic Higgs bundles and
irreducible
representations of the fundamental group of punctured Riemann surfaces
established by Simpson. We
generalize a result of Hitchin, identifying those parabolic
Higgs bundles that correspond to Fuchsian representations.
We also study the Higgs bundles that
give representations whose image is contained, after conjugation,
in SL($k,\Bbb R$). We compute the real dimension of one of the
components of this space of representations,
which in the absence of punctures is the generalized Teichm\"uller
space introduced by Hitchin,
and which in the case of $k=2$ is the Teichm\"uller space of the
Riemann surface.
\end{abstract}

\section{Introduction}
In the well-known paper \cite{H1}, Hitchin introduced Higgs
bundles, and established a one-to-one correspondence between
equivalence classes of irreducible GL($2,\Bbb{C}$) representations of
the fundamental group of a compact Riemann surface and
isomorphism classes of rank two stable Higgs of degree zero.
%Among other things in \cite{H2}, he identified the Higgs bundles
%corresponding to the Fuchsian representations.
In \cite{S2}, Simpson defined parabolic Higgs bundles, which
generalized Hitchin's correspondence to the case of
open Riemann surfaces (see also \cite{S1}).
%between Higgs
%bundles and representations to the case of an open Riemann
%surface, by introducing the notion of parabolic Higgs bundles.
Here, by an open Riemann surface we mean the complement of
finitely many points in a compact surface. More
precisely, Simpson identified what he calls filtered local
systems with parabolic Higgs bundles.

In \cite{H2}, Hitchin identified the Higgs bundles
corresponding to the Fuchsian representations.
Our main aim here is to generalize his results to the case of open
Riemann surfaces.
%identify the parabolic Higgs bundles
%corresponding to the Fuchsian representations of open Riemann
%surfaces.

Before giving more details, we describe the result of Hitchin on
Fuchsian representations.
Let $X$ be a compact Riemann surface of genus $g \geq 2$, and let
$L$ be a line bundle on $X$ such that $L^2 = K_X$, that is $L$
is a square root of the canonical bundle of $X$. Define $$ E
\, := \, L^* \oplus L$$ which is a rank $2$ bundle on $X$. For
$a\in H^0(X,K^2)$, let $${\theta}(a) \, :=\, \left(\begin{matrix}0 &
1\\ a & 0\\ \end{matrix}\right) \, \in \,H^0\big({\bar X},
End(E)\otimes K\big)$$ be the Higgs field. Hitchin proved that
the conjugacy classes of Fuchsian representations of $\pi_1(X)$
(homomorphisms of ${\pi}_1(X)$ into PSL$(2,\Bbb R$) such that the
quotient of the action on the upper half plane is a compact
Riemann surface of genus $g$) correspond to the Higgs bundles of
the form $(E,{\theta}(a))$ defined above. Moreover, the Higgs
bundle $(E,{\theta}(0))$ corresponds to the Fuchsian representation
for the Riemann surface $X$ itself.

%Now we will describe the parabolic Higgs bundles corresponding
%to the Fuchsian representations of the fundamental group of an
%open Riemann surface $X = {\bar X} -\{p_1,\ldots,p_n\}$, where $\bar X$
%is a compact Riemann surface of genus $g$, and $2g-2 +n >0$.
Consider now an open Riemann surface $X = {\bar X} - \{ p_1, \ldots, p_n
\}$, where $\bar X$ is a compact surface of genus $g$, and $p_1, \ldots,
p_n$ are $n>0$ distinct points of $\bar X$. Let $D$ denote the divisor
given by these points, \ie $D=\{p_1,\dots,p_n\}$.
We will further assume that
$2g-2+n>0$, which is equivalent to the condition that the
universal covering space of $X$ is (conformally equivalent to) the upper
half plane.
Consider the bundle $$E\,\, :=\,\, (L\otimes {\cal O}_X(D))^*\,
\oplus \,L $$ and give parabolic weight $1/2$ to the fiber
$E_{p_i}$, $1\leq i \leq n$. For $a\in H^0({\bar X},K^2\otimes
{\cal O}_{\bar X}(D))$ let $${\theta}(a) \, :=\, \left(\begin{matrix}0 & 1\\ a
& 0\\ \end{matrix}\right) \, \in \,H^0\big({\bar X}, End(E)\otimes
K\otimes {\cal O}_X(D)\big) $$ be a parabolic Higgs field on the
parabolic bundle $E$.

We prove that under the identification between filtered local
systems and parabolic Higgs bundles, Fuchsian representations of
$n$-punctured Riemann surfaces are in one-one correspondence
with parabolic Higgs bundles of the type $(E,{\theta}(a))$ defined above.
Moreover, the parabolic Higgs bundle $(E,{\theta}(0))$ corresponds to the
Fuchsian representation of the punctured surface $X$ itself.
Thus this is a direct generalization of the result of Hitchin on
Fuchsian representations of compact Riemann surfaces to the
punctured case.

In section $3$, we generalize the above results to the case of
representations of the fundamental group of the surface $X$ into
PSL($k,\Bbb R$), for $k>2$. More precisely, we consider a parabolic
bundle $W_k$,
obtained by tensoring the $(k-1)$th symmetric product of the bundle
$E$ defined above with an appropriate power of ${\cal O}_{\bar X}(D)$.
The Higgs fields we consider are generalizations of the $2$-dimensional
case, namely they are of the form
$$\theta(a_2,\ldots,a_{k-1}):=\left(\begin{matrix}0 & 1 & \cdots & 0 \\
0 & 0 & 1 & \vdots \\
\vdots &   &   & 0 \\
a_k & \cdots & a_2 & 0 \end{matrix}\right) ,$$
%where the $a_j$
%are $j$th differentials on $\bar X$, with at most poles
%of order $j-1$ at the points of $D$. Or in other words, the elements
$a_j$ is a section of the line bundle
$K^j \otimes ({\cal O}_{{\bar X}(D)})^{j-1}$.
As in section $2$, we have that
the pair $(W_k,\theta(a_2, \ldots, a_k))$ is a stable parabolic bundle
of parabolic degree $0$. It is not difficult to see that the parabolic
dual of $W_k$ is naturally isomorphic to
the parabolic bundle $W_k$ itself.  This
implies that the holonomy of the flat connection corresponding to these
bundles is contained in a split real form of SL($k,\Bbb C$),
which is isomorphic to SL($k,\Bbb R$).
%itself. This implies that we have an embedding of
%$\pi_1(X)$ into SL($k,\Bbb R$), by means of flat connections on $W_k$.
We prove that one of the components of the space of
representations of $\pi_1(X)$ into SL$(k,\Bbb R)$, with fixed
conjugacy class of monodromy
around the punctures, has real dimension equal to
$2(k^2-1)(g-1)+k(k-1)n$. Observe that for $k=2$, this dimension is
$2(3g-3+n)$, which is precisely the real dimension of ${\cal T}_g^n$,
the Teichm\"uller space of compact surfaces of genus $g$ with $n$
punctures. It is therefore natural, following \cite{H2}, to call the
above component the {\it Teichm\"uller component} of the corresponding
space of representations. Further study of this set is worthwhile.

\section{Higgs bundles for Fuchsian representations}
Let $\bar X$ be a compact Riemann surface of genus $g$, and let $$D \,\,
:= \,\, \{p_1,p_2,\ldots ,p_n\}$$ be $n$ distinct points on $\bar X$.
Define $X := {\bar X} - D$ to be the punctured Riemann surface given by
the complement of the divisor $D$.
We will assume that $2g-2+n > 0$, that is, the surface $X$ supports a
metric of constant curvature $(-4)$.

The degree of the holomorphic cotangent bundle $K$, of $\bar X$ is
$2g-2$. Therefore, there is a line bundle $L$ on $\bar X$ such that
$L^2 = K$. Fix such a line bundle $L$. Note that any two of the
$4^g$ possible choices of $L$ differ by a line bundle of order
$2$.
%Note that any two choices of $L$ differ by a line bundle of order 2.
%(There are $4^g$ of them.)

%Let us recall the definition of the parabolic Higgs bundles introduced
%in \cite{S2}. Using $L$ we will construct a parabolic Higgs
%bundle on $\bar X$.

Using $L$ we will construct a parabolic Higgs bundle on $\bar X$, as
follows.
Let $\xi = {\cal O}_{\bar X}(D)$ denote the line bundle on $\bar X$ given by
the divisor $D$. Define
\begin{equation} E\,\, :=\,\, (L\otimes \xi)^*\,
\oplus \,L \label{eq:bundle}\end{equation}
to be the rank 2 vector bundle on
$\bar X$. To define a parabolic structure on $E$
(we will follow the definition of parabolic Higgs bundle in \cite{S2}),
on each point
$p_i\in D$, $1\leq i \leq n$, we consider the trivial flag
$$ E_{p_i} \, \supset \, 0,$$
and give parabolic weight $1/2$ to
$E_{p_i}$. This gives a parabolic structure on $E$.

Note that
\begin{equation} Hom(L , L^{*}\otimes {\xi}^*)\otimes K\otimes\xi
\, =\, {\cal O} \, \subset \, End(E)\otimes K\otimes \xi
\label{eq:hom}\end{equation}
Let $1$ denote the section of $\cal O$ given by the
constant function $1$. So from (\ref{eq:hom}) we have
\begin{equation} \theta \, :=\,
\left(\begin{matrix}0 & 1\\ 0 & 0\\ \end{matrix}\right) \, \in
\,H^0\big({\bar X}, End(E)\otimes K\otimes\xi\big) \label{eq:theta}
\end{equation}

\begin{lemma} The parabolic Higgs bundle
$(E, \theta)$ is a parabolic stable Higgs bundle of parabolic
degree zero.\label{lemma:stable}\end{lemma}

\begin{pf} From the definition of parabolic degree
(see \cite[definition 1.11]{MS} or \cite{S2}) we immediately
conclude that the parabolic degree of $E$ is zero.

To see that $(E, \theta)$ is stable, first note that there is only
one sub-bundle of $E$ which is invariant under $\theta$, namely the summand
$(L\otimes \xi)^*$ in (2.1). (A sub-bundle $F \subset E$ is
called {\it invariant} under $\theta$ if $\theta (F) \subset F\otimes
K\otimes \xi$.) The degree of $(L\otimes \xi)^*$ is $1-g -n$. So
the parabolic degree of $(L\otimes \xi)^*$, for
the induced parabolic structure, is $1-g-n/2$.

But, from our assumption that $2g-2+n > 0$ we have $1-g-n/2 <
0$. So $(E,\theta)$ is stable.
\end{pf}

{}From the proof of the \lemref{lemma:stable} it follows that $(E,\theta)$
constructed above is stable if and only if $2g-2+n > 0$. We will
show later that this corresponds to the fact that $X$ admits a
complete metric of constant negative curvature if and only if
$2g-2+n > 0$.

{}From the main theorem of \cite[pg. 755]{S2} we know that there
is a tame harmonic metric on the bundle $E$. (See the Synopsis
of that paper for the definition of tame harmonic metric.)

It is well-known that there is an unique complete K\"ahler metric
on $X$, known as the {\it Poincar\'e} {\it metric},
such that its curvature is $(-4)$.

Both the bundles $L$ and $(L\otimes \xi)^*$ are equipped with
metrics induced by the tame harmonic metric on $E$. So $$
Hom(L,(L\otimes \xi)^*) \, =\, L^2\otimes{\xi}^* \, =\, T\otimes
\xi$$ is equipped with a metric. The restriction to $X$ of the line
bundle $\xi$, and hence ${\xi}^*$, on $\bar X$
has a canonical trivialization. Therefore we have a hermitian
metric on $T_X$ the tangent bundle of $X$. We will denote
this hermitian metric on $T_X$ by $H$. Note that $H$ is
singular at $D$, \ie does not induce a hermitian metric on
$T_{\bar X}$.

\begin{lemma} The hermitian metric $H$ on
the holomorphic tangent bundle on $X$ obtained above is the
Poincar\'e metric.
\label{lemma:metric}\end{lemma}

\begin{pf} We recall the Hermitian-Yang-Mills
equation which gives the harmonic metric on $E$ \cite{S2}. This
equation was first introduced in \cite{H1}.

Let $\nabla$ denote the holomorphic hermitian connection on the
restriction of $E$ to $X$ for the harmonic metric. Then the
Hermitian-Yang-Mills equation of the curvature of $\nabla$ is the following:
\begin{equation} K(\nabla) \, :=\, {\nabla}^2 \, = \, -\, [\theta,
{\theta}^*]\label{eq:yang}\end{equation}
If the decomposition (\ref{eq:bundle}) is orthogonal
with respect to the metric, then $[\theta, {\theta}^*]$ is a 2-form
with values in the diagonal endomorphisms of $E$ (diagonal for
the decomposition (\ref{eq:bundle})). Using this, the equation
(\ref{eq:yang})
reduces to the following equation on $X$
\begin{equation} F_H \,=\, -2{\bar H},
\label{eq:curv}\end{equation} where $H$ is a hermitian metric on $T_X$ and
$\bar H$ is the $(1,1)$-form on $X$ given by $H$. Observe that
$\bar H$ also denotes the K\"ahler $2$-form for the metric $H$.

A metric $H'$ on $T_X$ induces a metric on $L$. Since the bundle
$\xi$ has a natural trivialization over $X$, the metric $H'$
also induces a metric on $(L\otimes\xi)^*$, and therefore also
on $E$. If $H'$ satisfies the equation (\ref{eq:curv}) then the
metric on $E$ obtained this way satisfies (\ref{eq:yang}). Now from the
uniqueness of the solution of (\ref{eq:yang}) (\cite{S2}), we have that
such metric is obtained from the solution of (\ref{eq:curv}) in the above
fashion.
%the solution of (\ref{eq:yang}) is obtained from the solution of
%\ref{eq:curv} in the above fashion.

{}From the computation in Example (1.5) of \cite[pg. 66]{H1}, we
conclude that the K\"ahler metric $H$ on $X$ has Gaussian
curvature $(-4)$.

So in order to complete the proof of the lemma we must show that
the K\"ahler metric on $X$ is complete.

Recall the asymptotic behavior of the harmonic metric near the
punctures given in Section 7 of \cite{S2}. First of all, observe
that the fiber of $K\otimes \xi$ at any
$p_i \in D$ is canonically isomorphic to $\Bbb C$. So the fiber
$(End(E)\otimes K\otimes \xi)_{p_i}$ is $End(E_{p_i})$. The
evaluation of the section ${\theta}$ at $p_i$ as an element of
$End(E_{p_i})$ is defined to be the residue of $\theta$ at $p_i$.
For the Higgs field $\theta$, we have that the residue
at each $p_i$ is $$N \, := \, \left(\begin{matrix}0 & 1\\ 0 &
0\\ \end{matrix}\right).$$
In  \cite[pg. 755]{S2}, Simpson studies parabolic Higgs
bundles with residue $N$ as above. Consider the displayed
equation in page 758 of \cite{S2}, which describes the asymptotic
behavior of the harmonic metric. Using the fact that the
parabolic weight of $E_{p_i}$ is $1/2$ we conclude that for the
metric on $L$ induced by the tame harmonic metric on
$E$, both $a_i$ and $n_i$ in the equation in page 758 of
\cite{S2} are $1/2$. (We also use the fact that, in the notation
of \cite[pg. 755]{S2}, $L \subset W_1$ and $L$ is not contained in
$W_0$.) In other words, in a suitable trivialization of $L$
on an open set containing a puncture $p_i \in D$, and with
holomorphic coordinate $z$ around $p_i$, the hermitian metric on
$L$ obtained by restricting the harmonic metric on $E$ is
$$r^{1/2}|{\mathrm{log}}(r)|^{1/2},$$ where $r = |z|$.

Similarly, for $(L\otimes\xi)^*$, the $a_i$ and $n_i$ in
the equation \cite[pg. 758]{S2} are $1/2$ and $-1/2$
respectively.

So the metric on $Hom(L,(L\otimes\xi)^*)$ is $({\mathrm{
log}}|(r)|)^{-1}$. Recall the earlier remark that ${\xi}^*$ has a
natural trivialization on $X$. The section of ${\xi}^*$ on
$X$ has a pole of order $1$ at the points of $D$,
when it is considered as a
meromorphic section of ${\xi}^*$ on $\bar X$. This implies that the
hermitian metric on $T = L^{-2}$ is
\begin{equation} r^{-1}|{\mathrm{log}}(r)|^{-1}.
\label{eq:metric}\end{equation}
But this is the expression the Poincar\'e metric
of the punctured disk in $\Bbb C$. This proves that the K\"ahler
metric on $X$ induced by $H$ is indeed complete. This completes
the proof of the lemma.
\end{pf}

{}From the decomposition (\ref{eq:bundle}) it follows that
\begin{equation} Hom(L^{*}\otimes
{\xi}^* ,L)\otimes K\otimes\xi \, =\, K^2\otimes {\xi}^2 \,
\subset \, End(E)\otimes K\otimes \xi
\label{eq:hom2}\end{equation}

Note that the bundle $\xi$ has a natural section, which we will denote
by $1_{\xi}$. We may imbedd $H^0({\bar X} ,K^2\otimes \xi)$ into
$H^0({\bar X} ,K^2\otimes {\xi}^2)$ by the homomorphism $s \longmapsto
s\otimes 1_{\xi}$. So using (\ref{eq:hom2}) we have a natural homomorphism
\begin{equation} \rho \, :\, H^0({\bar X}, K^2\otimes \xi) \,
\longrightarrow \, H^0({\bar X} ,End(E)\otimes K\otimes \xi)
\label{eq:rho} \end{equation}
Note that the image
of $\rho$ is contained in the image of the inclusion $$H^0({\bar X},
End(E)\otimes K) \, \longrightarrow \,H^0({\bar X},
End(E)\otimes K\otimes \xi)$$ With a slight abuse of notation,
for any $a\in H^0({\bar X},K^2\otimes \xi)$, the corresponding element in
$H^0({\bar X} ,End(E)\otimes K)$ will also be denoted by $\rho (a)$.

The following theorem is a generalization of theorem (11.2) of
\cite{H1} to the case of open Riemann surfaces.
\begin{thm} For any $a \in
H^0({\bar X} ,K^2\otimes \xi)$, the Higgs structure $${\theta}_a \,
:=\,{\theta}
+ {\rho}(a) \, =\, \left(\begin{matrix}0 & 1\\ 0 & 0\\
\end{matrix}\right) + {\rho}(a)$$ on the parabolic bundles $E$
(defined in (\ref{eq:bundle})) makes $(E,{\theta}_a)$ a parabolic stable Higgs
bundle of parabolic degree zero.

Let $H_a$ denote the harmonic metric (given by the main theorem of
\cite{S2}) on the restriction of $E$ to $X$, and let $h$ denote
the K\"ahler metric on $X$ induced by the tame harmonic metric
${H}_{a}$ as in \lemref{lemma:metric}. Then the following holds :

\begin{enumerate}
\item{} The section of the $2$-nd symmetric power of the complex
tangent bundle $$h_a \, :=\, a+h+{\bar a}+a{\bar a}/h \, \in
\, {\Omega}^0(X, S^2T^*\otimes\Bbb C)$$ is a Riemannian metric on $X$.

\item{} The metric $h_a$ is a complete Riemannian metric of
constant Gaussian curvature $(-4)$. The Riemann surface
structure on $X$ given by metric $h_a$ is a Riemann surface with
punctures, \ie there are no holes. (A Riemann surface with a hole
is a complement of a disk in a compact Riemann surface.)

\item{} Associating to $a \in H^0({\bar X} ,K^2\otimes \xi)$ the complex
structure on the $C^{\infty}$ surface $X$ given by the metric
$h_a$, the map obtained from $H^0({\bar X} ,K^2\otimes \xi)$ to the
Teichm\"uller space ${\cal T}^n_g$ of surfaces of genus $g$ and $n$
punctures
%(there are no holes, by statement $2$ above),
is a bijection.
\end{enumerate}
\label{thm:main}\end{thm}

\begin{pf} To prove that $(E,{\theta}_a)$ is stable we
use a trick of \cite{H2}. For $\mu >0$, define an
automorphism of $E$ by $$T\, := \, \left(
\begin{matrix}1 & 0\\ 0 & \mu\\ \end{matrix}\right) .$$
The parabolic Higgs bundle $(E,{\theta}_a)$ is
isomorphic to $(E,T^{-1}\circ {\theta}_a\circ T)$, and hence
$(E,T^{-1}\circ {\theta}_a\circ T)$ is parabolic stable if and
only if $(E,{\theta}_a)$ is so. Since $\mu \neq 0$, we have
$(E,T^{-1}\circ {\theta}_a\circ T)$ is parabolic stable if and only
if $(E,\frac{1}{\mu} T^{-1}\circ {\theta}_a\circ T)$ is parabolic stable.
Now $$1/\mu T^{-1}\circ {\theta}_a\circ T\,=\,\left(\begin{matrix}
0 & 1 \\ 0 & 0\\ \end{matrix}\right) +{\rho}(a) /\mu \, = \,
{\theta}_{a/\mu}.$$ But from the openness of the stability
condition we have that since $(E, \theta)$ is stable [\lemref{lemma:stable}],
there is a non-empty open set $U$ in $H^0({\bar X} ,K^2\otimes \xi)$
containing the origin such that for any $a\in U$, the bundle
$(E,{\theta}_a)$ is parabolic stable. Taking $\mu$ to be sufficiently
large so that ${\theta}_{a/\mu} \in U$, we conclude that any
$(E,{\theta}_a)$ is parabolic stable.

The bundle $E$ is equipped with the harmonic metric $H_a$, and
$K$ has a metric induced by $h_a$. Using these metrics we
construct a hermitian metric on $End(E)\otimes K$. Since $\rho(a)
\in H^0({\bar X} ,End(E)\otimes K)$, we may take its pointwise norm.

To prove the statement ($1$) we first want to calculate the
behavior of $||\rho(a)||$ near the punctures. Since $\rho(a)
\in H^0({\bar X} ,End(E)\otimes K)$, we have $$
{\mathrm{residue}}\, ({\theta}_a)\, =\, {\mathrm{residue}}\,
({\theta}) \, =\, N.$$ So the
two hermitian metrics $H_0$ (corresponding to $a=0$) and $H_a$
on $E$ are mutually bounded, \ie $C_1.H_0 \leq H_a \leq C_2.H_0$
for some constants $C_1$ and $C_2$. (Recall that the metric in
\lemref{lemma:metric} was induced by $H_0$.) From this it is easy to check
that around any puncture $p_i$, the norm $||\rho(a)||$ is bounded
by $r|{\mathrm{log}}(r)|^{3/2}$. This implies that $||\rho(a)||$
converges to zero as we approach a puncture.

%Now we are in a position to imitate the argument in theorem
%($11.2$) of \cite{H1} to prove the statement (1).
Arguing as in
($11.2$) of \cite{H1}, if $h_a$ is not a metric then $$1\, - \,
||\rho(a)|| \, \leq\, 0$$ at some point $x\in X$. Since $||\rho(a)||$
converges to zero as we approach a puncture, the infimum of the
function $1 - ||\rho(a)||$ on $X$ must be attained somewhere, say
at $x_0\in X$.

Let $\Delta$ denote the Laplacian operator acting on smooth
functions on
$X$. Since the operator ${\cal L} := - \Delta -4.||\rho(a)||^2$ is
uniformly elliptic on $X$, we may apply \cite[Section VI.3.,
Proposition 3.3]{JT} for the operator ${\cal L}$ and the point
$x_0$.
%From [JT, Section VI.3., Proposition 3.3]
We conclude that either $1 - ||\rho(a)|| >0$ or $1 - ||\rho(a)||$
is a constant
function. This proves that $h_a$ is a Riemannian metric on $X$.

{}From the computation in the proof of Theorem (11.3)(ii) of
\cite[pg. 120]{H1}, we conclude that
$h_a$ is a metric of curvature $(-4)$.
%the curvature of $h_a$ is
%$(-4)$.

To complete the proof of the statement ($2$) we must show that
$h_a$ is complete and it has finite volume. (If the volume of
the Poincar\'e metric on a Riemann surface is finite then the
Riemann surface is a complement of finite number of points in a
compact Riemann surface. In particular, the Riemann surface can
not have any holes.)

The above established fact that the metrics $H_0$ and $H_a$
on $E$ are mutually bounded, together with \lemref{lemma:metric} imply
that the Riemannian metric $h_a$ and the Poincar\'e metric on
$X$ are mutually bounded. Since the Poincar\'e metric is
complete and of finite volume, the same must hold for $h_a$.

To prove the statement ($3$) we have to show that map from $H^0({\bar X}
,K^2\otimes \xi)$ to the Teichm\"uller space ${\cal T}^n_g$ obtained
in ($2$) is surjective. This will follow from Section $3$ where we
will prove that the image is both open and closed, and hence it
must be surjective as ${\cal T}^n_g$ is connected.

However we may also use the argument in \cite[Theorem
(11.2)(iii)]{H1} to prove statement (3). Let $h_0$ denote the
Poincar\'e metric on $X$. Indeed, to make the argument work all
we need to show is the following generalization of the
Eells-Sampson theorem to punctured Riemann surfaces:
given a complete Riemannian
metric $h$ of constant curvature $(-4)$ and finite volume on the
$C^{\infty}$ surface $X$, there is a unique diffeomorphism $f$,
of $X$ homotopic to the identity map, such that $f$ is a harmonic
map from $(X,h_0)$ to $(X,h)$.
%`(The generalization of the
%`Eells-Sampson theorem to punctured Riemann surfaces.)
This follows from the generalization of the theorem of Corlette,
\cite{C}, to the non-compact case as mentioned in \cite[pg. 754]{S2}.

Let $(V,\nabla)$ be the flat rank two bundle given by the Fuchsian
representation for the Riemann surface $(X,g)$. Let $H$ be the
harmonic metric on $V$ given by the main theorem of \cite{S2}
(pg. $755$) for the flat bundle $(V,\nabla)$ on the Riemann surface
$(X,h_0)$. In other words, $H$ gives a section, denoted by $s$,
of the associated bundle with fiber SL$(2,{\Bbb R})/$SO($2$) = $\Bbb H$,
where $\Bbb H$ is the upper half plane. This section $s$ gives the
harmonic map $f$ mentioned above. This completes the proof of
the theorem.
\end{pf}

The vector space $H^0({\bar X} ,K^2\otimes \xi)$ has a natural complex
structure. So does the Teichm\"uller space ${\cal T}_g^n$.
%The Teichm\"uller space ${\cal T}^n_g$ has a natural
%complex structure.
The identification of $H^0({\bar X} ,K^2\otimes \xi)$
with ${\cal T}^n_g$ given by Theorem 2.11 does not preserve the
complex structures. Indeed, ${\cal T}^n_g$ is known to be
biholomorphic to a bounded domain in ${\Bbb C}^{3g-3+n}$. Since any
bounded holomorphic function on an affine space must be
constant, the identification in Theorem 2.11 is never
holomorphic.

\noindent {\bf Remark}\, The parabolic dual of the
parabolic bundle $E$ is $E^*\otimes {\xi}^*$ with trivial
parabolic flag and parabolic weight $1/2$ at the parabolic
points $p_i$, $1\leq i \leq n$.  So the parabolic dual of $E$ is
$E$ itself. Any parabolic Higgs bundle $(E, {\theta}_a)$ (as in
\thmref{thm:main}) is naturally isomorphic to the parabolic Higgs
bundle $(E^*, {\theta}^*_a)$, where $E^*$ is the parabolic dual of
$E$. This implies that the holonomy of the flat connection on
$X$ corresponding to the Higgs bundle $(E,{\theta}_a)$ is contained
(after conjugation) in SL$(2,\Bbb R$). This of course is also
implied by \thmref{thm:main} since the image of a Fuchsian
representation is contained in PSL$(2,\Bbb R$).

\section{Higgs bundles for
SL(${\load{\normalsize}{\it}k},\Bbb R$) representations}
Recall the bundle $E$ of section 2, which was defined by $E = (L \otimes
\xi)^* \oplus L$, where $L$ is a (fixed) square root of the canonical
bundle $K$, and $\xi={\cal O}_{\bar X}(D)$. The ($k-1$)-th
symmetric product of ${\Bbb C}^2$ produces an embedding of
SL($2,\Bbb R$) into SL($k,\Bbb R$), via action on homogeneous polynomials of
degree $k$. Let $V_k$ denote the bundle given by the ($k-1$)-th
symmetric product of $E$, that is $V_k:=S^{k-1}(E)$.
At each point $p_i\in D$ we have the trivial flag
$(V_k)_{p_i}\supset 0$, $1\leq p_i \leq n$, with weight
equal to $\frac{k-1}{2}$. In order to construct a parabolic bundle, we
need to reduce the weight to a number in the interval $[0,1)$. We do
this by tensoring $V_k$ with $\xi^{m(k)}$, where $m(k)$ is equal to
$\frac{k}{2}-1$, if $k$ is even, or $\frac{k-1}{2}$, if $k$ is odd. We
will denote the bundle $V_k\otimes\xi^{m(k)}$ by $W_k$. At each point
$p_i\in D$,  we take the trivial flag $(W_k)_{p_i} \supset 0$
of $W_k$, with weight equal to
$\frac{1}{2}$, if $k$ is even, or $0$, if $k$ is odd.

Considering $1$ as the section of $\cal O$ given by the
constant function $1$, we can define
\begin{equation}
\theta(0,\ldots,0):=\left(\begin{matrix}0 & 1 & \cdots & 0 \\
0 & 0 & 1 & \vdots \\
\vdots &  &  & 1 \\
0 & 0 & \cdots & 0 \end{matrix}\right) ,\label{eqn:theta0}\end{equation}
which represents an element of $H^0\big( {\bar X}, End(W) \otimes K \otimes
\xi)$.

\begin{lemma}The bundle $(W_k,\theta(0,\ldots,0))$ is a parabolic stable
Higgs bundle of para-\newline bolic degree zero.\end{lemma}
\begin{pf}
If $k$ is even, we have that the parabolic degree of $W_k$ is
equal to $\frac{k(k+1)}{2}n + \frac{k}{2}(k+1)n=0$.
In the case of odd $k$,
it is easy to see that the degree (as a bundle)
of $W_k$ is $0$, and since the weight is equal to $0$,
we get that the parabolic degree of $W_k$ is zero.

The invariant proper sub-bundles of \ref{eqn:theta0} are
$L^{1-k}\otimes \xi^{-k/2},\ldots,
L^{1-k}\otimes \xi^{-k/2}\oplus\cdots\oplus
L^{(k/2)-1}\otimes\xi^{(k-4)/2}$, if
$k$ is even; or
$L^{1-k}\otimes \xi^{(1-k)/2},\ldots,
L^{1-k}\otimes \xi^{(1-k)/2}\oplus\cdots\oplus
L^{k-3}\otimes\xi^{(k-3)/2}$,
if $k$ is odd.
It is not difficult to see that all these sub-bundles have negative
parabolic degree.
\end{pf}

Using the natural section $1_\xi$ of $\xi$, we embed the spaces
$H^0({\bar X},K^j \otimes \xi^{j-1})$, $j=2,\ldots,k$, into $H^0({\bar
X},End(W_k)\otimes K \otimes \xi)$.
By an abuse of notation, if $a_j \in H^0({\bar X},K^j \otimes
\xi^{j-1})$, we understand the above embedding as producing an element
\begin{equation}
\theta(a_2,\ldots,a_{k-1}):=\left(\begin{matrix}0 & 1 &
\cdots & 0 \\
0 & 0 & 1 & \vdots \\
\vdots &   &  & 1 \\
a_k & \cdots & a_2 & 0 \end{matrix}\right)\label{eqn:thetaa}\end{equation}
of $H^0({\bar X},End(W_k) \otimes K \otimes \xi)$.
Now, by the arguments of Hitchin, based on the
openness of the stability of bundles, we get that the pair $(W_k,
\theta (a_2,\ldots,a_k))$ is a stable parabolic Higgs bundle of
parabolic degree $0$.
Using these special Higgs bundles, one can obtain some information about
the space of representations of the fundamental group of $X$ into
SL($k,\Bbb R$). More precisely, our result is as follows.
\begin{prop}
The space of representations of the fundamental group of $X$ in
SL($k,\Bbb R$), with fixed conjugacy class of monodromy around the punctures,
has a component of real dimension
$2(k^2-1)(g-1)+k(k-1)n$.
\end{prop}
\begin{pf}
By the work of Simpson \cite{S2} and Balaji Srinivasan
\cite{B}, we have a one-to-one continuous
correspondence
between the space $M$ of stable parabolic Higgs bundles of degree zero,
and the space of representations of the fundamental group of $X$ into
SL($k,\Bbb C$).
Consider the parabolic dual of $W_k$, which is constructed as follows.
First, take the dual bundle $W_k^*$ of $W_k$.
If $k$ is odd, since the weight of
the flag is $0$, we have that the parabolic dual of $W_k$ is $W_k^*$,
with trivial flag at the points $p_i\in D$, and weight equal to zero.
If $k$ is even, we have a weight of $-\frac{1}{2}$ associated to the
trivial flag of $W_k^*$.
%To bring this weight up to a number between $0$
%and $1$, we
Tensor $W_k^*$ with $\xi$ to obtain that the parabolic dual
of $W_k$ is  $W_k^* \otimes \xi$. So we always have that the parabolic
dual of the bundle $W_k$ is $W_k$ itself. This implies that the image of
the fundamental group under the representation induced by $(W_k,\theta)$
lies in SL($k,\Bbb R$).

Since $a_j$ is a section of  $K^j \otimes \xi^{j-1}$, we have that the residue
of the Higgs field is invariant, \ie
$${\mathrm{residue}}\,
(\theta(a_2,\ldots,a_{k-1})) =
{\mathrm{residue}}\, (\theta(0,\ldots,0)) =
\left(\begin{matrix}0 & 1 & \cdots & 0 \\
0 & 0 & 1 & \vdots \\
\vdots &  &  & 1 \\
0 & 0 & \cdots & 0 \end{matrix}\right) .$$ This implies that in the above
representation, the conjugacy class of the elements corresponding to
small loops around the punctures of $X$
is invariant. By the embedding of SL($2,\Bbb R$) into SL($k,\Bbb R$),
we have that this is the class of the element
\begin{equation}
\left(\begin{matrix}1 & 1 & \cdots & 0 \\
0 & 1 & 1 & \vdots \\
\vdots &  &  & 1 \\
0 & 0 & \cdots & 1 \end{matrix}\right) .\end{equation}

Using a bases, $\{p_1,\ldots,p_{k-1}\}$,
for the set of invariant polynomials of the
Lie algebra of SL($k,\Bbb C$), we can construct a continuous mapping
$p:M\rightarrow \bigoplus_{j=2}^{k}H^0(\bar{X},K^j\otimes\xi^j)$, given by
assigning to the Higgs field $(W_k,\Phi)$ the elements
$(p_1(\Phi),\ldots,p_{k-1}(\Phi))$. The Higgs fields of the form
(\ref{eqn:thetaa})
produce a section $s$ of $p$, defined over the closed subspace
$\bigoplus_{j=2}^{k}H^0(\bar{X},K^j\otimes\xi^{j-1})$.
Therefore, we have that the image of $s$ is closed.
One can easily compute that the dimension
(over $\Bbb R$) of the space of sections
$\bigoplus_{j=2}^{k}H^0(\bar{X},K^j\otimes\xi^{j-1})$ is equal to
$$\sum_{j=2}^k 2(2j-1)(g-1)+2\sum_{j=2}^k(j-1)n=2(k^2-1)(g-1)+k(k-1)n.$$
On the
other hand, the dimension of the space of representations of the
fundamental group of $X$ into SL($k,\Bbb R$), with the condition that the
monodromy around the punctures lies in the above conjugacy class, can be
computed as follows. The fundamental group of $X$ can be identified
with a group of M\"obius transformations (or elements of SL($2,\Bbb
R$)), generated by elements
$\{c_1,d_1,\ldots,c_g,d_g,e_1,\ldots,e_n\}$, and with one
relation of the form $\prod_{j=1}^g[c_j,d_j]\prod_{j=1}^n e_j=id$, where
$[c,d]=cdc^{-1}d^{-1}$ denotes the commutator of the elements $c$ and $d$.
In classical terms, the transformations $c_j$'s and $d_j$'s are hyperbolic,
that is conjugate to dilatations, while the $e_j$'s are parabolic, or
conjugate to translations. In terms of loops on $X$, we have that the
$c_j$'s and $d_j$'s can be identified with paths around the handles of
$X$, while the $e_j$'s are simple loops around the punctures.
The image of the elements
$c_j$ and $d_j$ depends on dim(SL($k,\Bbb R$))=$k^2-1$ parameters.
%Let $U$ denote the set of regular unipotent elements of SL($k,\Bbb R$).
%An
%elements of $U$ depends on $\frac{k(k-1)}{2}$ real parameters.
In order to compute the number of parameters of the elements $e_j$,
first observe that these transformations belong to the conjugacy classes
of elements of $U=\{$regular unipotent elements of PSL($2, \Bbb
R$)$\}$.
Any matrix of SL($k,\Bbb R$) can be written as $ldu$, where $l$ is
unipotent and lower triangular, $d$ is diagonal, and $u$ is unipotent
upper triangular. We therefore have $(ldu)U(ldu)^{-1}=lUl^{-1}$. So the
conjugacy class of $U$ depends on $k(k-1)$ parameters.
Therefore, we have that the real dimensions of
$\bigoplus_{j=2}^{k}H^0(\bar{X},K^j\otimes\xi^{j-1})$ and the space of
representations of $\pi_1(X)$, with fixed conjugacy class for the
monodromy elements around the punctures, are equal.
Standard arguments using the invariance of domain theorem complete the
proof.\end{pf}
\ifx\undefined\bysame
\newcommand{\bysame}{\leavevmode\hbox to3em{\hrulefill}\,}
\fi


\begin{thebibliography}{1}

\bibitem{C}
K.~Corlette, {\em Flat {G}-bundles with canonical metrics}, J. Differential
  Geom. {\bf 28} (1988), 361--382.

\bibitem{H1}
N.~Hitchin, {\em The self-duality equation on a {R}iemann surface}, Proc.
  London Math. Soc. {\bf 55} (1987), 59--126.

\bibitem{H2}
\bysame, {\em Lie groups and {T}eichm{\"{u}}ller spaces}, Topology {\bf 31}
  (1992), 449--473.

\bibitem{JT}
A.~Jaffe and C.~Taubes, {\em Vortices and {M}onopoles}, Progress in {P}hysics,
  vol.~2, Birkh{\"{a}}user, Boston, {B}asel and {S}tuttgart, 1980.

\bibitem{MS}
V.~Mehta and C.~Seshadri, {\em Moduli of vector bundles on curves with
  parabolic structures}, Math. Ann. {\bf 248} (1980), 205--239.

\bibitem{S1}
C.~Simpson, {\em Constructing variations of {H}odge structure using
  {Y}ang-{M}ill theorey and application to uniformization}, J. Amer. Math. Soc.
  {\bf 1} (1988), 867--918.

\bibitem{S2}
\bysame, {\em Harmonic bundles on noncompact curves}, J. Amer. Math. Soc. {\bf
  3} (1990), 713--770.

\bibitem{B}
Balaji Srinivasan, Ph.D. thesis, University of {C}hicago, 1994.

\end{thebibliography}
\end{document}